\documentclass[english,journal]{IEEEtran}
\usepackage[T1]{fontenc}
\usepackage[latin9]{inputenc}
\usepackage{color}
\usepackage{amsthm}
\usepackage{amsmath}
\usepackage{amssymb}
\usepackage{esint}

\makeatletter
  \theoremstyle{definition}
  \newtheorem{defn}{\protect\definitionname}
  \theoremstyle{plain}
  \newtheorem{lem}{\protect\lemmaname}
  \theoremstyle{remark}
  \newtheorem{rem}{\protect\remarkname}
  \theoremstyle{plain}
  \newtheorem{cor}{\protect\corollaryname}

\pdfoutput=1
\usepackage{cite}

\makeatother

\usepackage{babel}
\providecommand{\corollaryname}{Corollary}
\providecommand{\definitionname}{Definition}
\providecommand{\lemmaname}{Lemma}
\providecommand{\remarkname}{Remark}

\begin{document}

\title{A Corollary for Nonsmooth Systems}

\author{N. Fischer, R. Kamalapurkar, W. E. Dixon}
\maketitle
\begin{abstract}
In this note, two generalized corollaries to the LaSalle-Yoshizawa
Theorem are presented for nonautonomous systems described by nonlinear
differential equations with discontinuous right-hand sides. Lyapunov-based
analysis methods are developed using differential inclusions to achieve
asymptotic convergence when the candidate Lyapunov derivative is upper
bounded by a negative semi-definite function.
\end{abstract}

\section{Introduction}

For continuous systems, stability techniques such as the LaSalle-Yoshizawa
Theorem provide a convenient analysis tool when the candidate Lyapunov
function derivative is upper bounded by a negative semi-definite function.
However, adapting the LaSalle-Yoshizawa Theorem to systems where the
time derivative of the system states are not locally Lipschitz remains
an open problem. The concept of utilizing the LaSalle-Yoshizawa Theorem
for nonsmooth systems was introduced in \cite{Nakakuki2008} as a
remark, but no formal proof was provided. 

In this note, we consider Filippov solutions for nonautonomous nonlinear
systems with right-hand side discontinuities%
\footnote{Throughout the subsequent presentation, a discontinuous right-hand
side will refer to being discontinuous in $x$, and continuous in
$t$.%
} utilizing Lipschitz continuous and regular Lyapunov functions whose
time derivatives (in the sense of Filippov) can be upper bounded by
negative semi-definite functions.

\section{Preliminaries}

Consider the system
\begin{equation}
\dot{x}=f\left(x,t\right)\label{eq: dynamics}
\end{equation}
where $x\left(t\right)\in D\subset\mathbb{R}^{n}$ denotes the state
vector, and $f:\:\mathcal{D}\times\left[0,\infty\right)\rightarrow\mathbb{R}^{n}$
is a Lebesgue measurable and essentially locally bounded \cite{Necula.Popescu.ea2010}
function. As is standard in literature\cite{Khalil1996}, existence
and uniqueness of the continuous solution $x\left(t\right)$ are provided
under the condition that the function $f$ is Lipschitz continuous.
However, if $f$ contains a discontinuity at any point in $\mathcal{D}$,
then a solution to (\ref{eq: dynamics}) may not exist in the classical
sense. Thus, it is necessary to redefine the concept of a solution.
Utilizing differential inclusions, the value of a generalized solution
(e.g., Filippov \cite{Filippov1988} or Krasovskii \cite{Krasovskii1963}
solutions) at a certain point can be found by interpreting the behavior
of its derivative at nearby points. Generalized solutions will be
close to the trajectories of the actual system since they are a limit
of solutions of ordinary differential equations with a continuous
right-hand side \cite{Guo2009}. While there exists a Filippov solution
for any arbitrary initial condition $x\left(t_{0}\right)\in\mathcal{D}$,
the solution is generally not unique \cite{Filippov1988,Aubin1984}.
\begin{defn}
\textbf{(Filippov Solution)} \label{def: filippov}\cite{Filippov1988}
A function $x\left(t\right)$ is called a solution of (\ref{eq: dynamics})
on the interval $\left[0,\infty\right)$ if $x\left(t\right)$ is
absolutely continuous and for almost all $t\in\left[0,\infty\right)$,
\[
\dot{x}\in K\left[f\right]\left(x,t\right)
\]
where
\begin{equation}
K\left[f\right]\left(x,t\right)\triangleq\bigcap_{\delta>0}\bigcap_{\mu N=0}\overline{co}f\left(B\left(x,\delta\right)\backslash N,t\right),\label{eq: K-def}
\end{equation}
$\underset{\mu N=0}{\bigcap}$ denotes the intersection over sets
$N$ of Lebesgue measure zero, $\overline{co}$ denotes convex closure,
and $B\left(x,\delta\right)=\left\{ \upsilon\in\mathbb{R}^{n}|\:\left\Vert x-\upsilon\right\Vert <\delta\right\} $.
\end{defn}
To facilitate the main results, three definitions are provided.
\begin{defn}
\textbf{(Directional Derivative) }\cite{Kaplan1991} Given a function
$f:\mathbb{R}^{m}\rightarrow\mathbb{R}^{n}$, the right directional
derivative of $f$ at $x\in\mathbb{R}^{m}$ in the direction of $v\in\mathbb{R}^{m}$
is defined as
\[
f'\left(x,v\right)=\lim\limits _{t\rightarrow0^{+}}\frac{f\left(x+tv\right)-f\left(x\right)}{t}.
\]
Additionally, the generalized directional derivative of $f$ at $x$
in the direction of $v$ is defined as
\[
f^{o}\left(x,v\right)=\lim\limits _{y\rightarrow x}\sup\limits _{t\rightarrow0^{+}}\frac{f\left(y+tv\right)-f\left(y\right)}{t}.
\]

\end{defn}

\begin{defn}
\textbf{(Regular Function)} A function $f:\mathbb{R}^{m}\rightarrow\mathbb{R}^{n}$
is said to be regular at $x\in\mathbb{R}^{m}$ if for all $v\in\mathbb{R}^{m}$,
the right directional derivative of $f$ at $x$ in the direction
of $v$ exists and $f'\left(x,v\right)=f^{o}\left(x,v\right)$.%
\footnote{Note that any $\mathcal{C}^{1}$ continuous function is regular and
the sum of regular functions is regular \cite{Clarke.Ledyaev.ea1998}.%
}
\end{defn}

\begin{defn}
\textbf{(Clarke's Generalized Gradient) }\cite{clarke1983} For a
function $V:\:\mathbb{R}^{n}\times\left[0,\infty\right)\rightarrow\mathbb{R}$
that is locally Lipschitz in $\left(x,t\right)$, define the generalized
gradient of $V$ at $\left(x,t\right)$ by
\[
\partial V\!\!\left(x,t\right)\!=\!\overline{co}\left\{ \lim\!\nabla V\!\!\left(x,t\right)\:\!|\:\!\left(x_{i},t_{i}\right)\!\rightarrow\!\left(x,t\right)\!,\left(x_{i},t_{i}\right)\notin\Omega_{V}\right\} 
\]
where $\Omega_{V}$ is the set of measure zero where the gradient
of $V$ is not defined.
\end{defn}
The following lemma provides a method for computing the time derivative
of a regular function $V\left(x,t\right)$ using Clarke's generalized
gradient \cite{clarke1983} and $K\left[f\right]\left(x,t\right)$,
from (\ref{eq: K-def}), along the solution trajectories of (\ref{eq: dynamics}).
\begin{lem}
\label{lem: chain rule}\textbf{\emph{(Chain Rule)}} \cite{Paden1987,Shevitz1994}
Let $x\left(t\right)$ be a Filippov solution of (\ref{eq: dynamics})
and $V:\:\mathcal{D}\times\left[0,\infty\right)\rightarrow\mathbb{R}$
be a locally Lipschitz, regular function. Then $V\left(x\left(t\right),t\right)$
is absolutely continuous, $\frac{d}{dt}V\left(x\left(t\right),t\right)$
exists almost everywhere (a.e.), i.e., for almost all $t\in\left[0,\infty\right)$,
and $\dot{V}\!\!\left(x\!\left(t\right)\!,t\right)\overset{a.e.}{\in}\dot{\tilde{V}}\!\!\left(x\!\left(t\right)\!,t\right)$,
where
\[
\dot{\tilde{V}}\left(x,t\right)\triangleq\bigcap_{\xi\in\partial V\left(x,t\right)}\xi^{T}\left(\begin{array}{c}
K\left[f\right]\left(x,t\right)\\
1
\end{array}\right).
\]
\end{lem}
\begin{rem}
Throughout the subsequent discussion, for brevity of notation, let
a.e. refer to almost all $t\in\left[0,\infty\right)$.
\end{rem}

\section{Main Result}

For the system described in (\ref{eq: dynamics}) with a continuous
right-hand side, existing Lyapunov theory can be used to examine the
stability of the closed-loop system using continuous techniques such
as those described in \cite{Khalil1996}. However, these theorems
must be altered for the set-valued map $\dot{\tilde{V}}\left(x\left(t\right),t\right)$
for systems with right-hand sides which are not Lipschitz continuous
\cite{Shevitz1994,Cheng2009,Guo2009}. Lyapunov analysis for nonsmooth
systems is analogous to the analysis used for continuous systems.
The differences are that differential equations are replaced with
inclusions, gradients are replaced with generalized gradients, and
points are replaced with sets throughout the analysis. The following
presentation and subsequent proofs demonstrate how the LaSalle-Yoshizawa
Theorem can be adapted for such systems.

The following auxiliary lemma from \cite{Paden1987} and Barbalat's
Lemma are provided to facilitate the proofs of the nonsmooth LaSalle-Yoshizawa
Corollaries.
\begin{lem}
\cite{Paden1987}\label{lem: v stable} Let $x\left(t\right)$ be
any Filippov solution to the system in (\ref{eq: dynamics}) and $V:\:\mathcal{D}\times\left[0,\infty\right)\rightarrow\mathbb{R}$
be a locally Lipschitz, regular function. If $\dot{V}\left(x\left(t\right),t\right)\overset{a.e.}{\leq}0$,
then $V\left(x\left(t\right),t\right)\leq V\left(x\left(t_{0}\right),t_{0}\right)\:\forall t>t_{0}$.\end{lem}
\begin{IEEEproof}
For the sake of contradiction, let there exist some $t>t_{0}$ such
that $V\left(x\left(t\right),t\right)>V\left(x\left(t\right),t_{0}\right)$.
Then,
\[
\int_{t_{0}}^{t}\dot{V}\left(x\left(\sigma\right),\sigma\right)d\sigma=V\left(x\left(t\right),t\right)-V\left(x\left(t\right),t_{0}\right)>0.
\]
It follows that $\dot{V}\left(x\left(t\right),t\right)>0$ on a set
of positive Lebesgue measure, which contradicts that $\dot{V}\left(x\left(t\right),t\right)\leq0$,
a.e. 
\end{IEEEproof}

\begin{lem}
\label{lem: barbalats}\textbf{\emph{(Barbalat's Lemma) }}\cite{Khalil1996}
Let $\phi:\:\mathbb{R}\rightarrow\mathbb{R}$ be a uniformly continuous
function on $\left[0,\infty\right)$. Suppose that $\lim\limits _{t\rightarrow\infty}\int_{0}^{t}\phi\left(\tau\right)d\tau$
exists and is finite. Then,
\[
\phi\left(t\right)\rightarrow0\:\:\: as\:\:\: t\rightarrow\infty.
\]

\end{lem}
Based on Lemmas \ref{lem: v stable} and \ref{lem: barbalats}, a
nonsmooth corollary to the LaSalle-Yoshizawa Theorem (c.f., \cite[Theorem 8.4]{Khalil2002}
and \cite[Theorem A.8]{Krstic1995}) is provided in Corollary \ref{thm: main_result}.
\begin{cor}
\label{thm: main_result}For the system given in (\ref{eq: dynamics}),
let $\mathcal{D}\subset\mathbb{R}^{n}$ be a domain containing $x=0$
and suppose $f$ is Lebesgue measurable and essentially locally bounded
on $\mathcal{D}\times\left[0,\infty\right)$. Furthermore, suppose
$f\left(0,t\right)$ is uniformly bounded for all $t\geq0$. Let $V:\mathcal{D}\times\left[0,\infty\right)\rightarrow\mathbb{R}$
be continuously differentiable in x, locally Lipschitz in t, and regular
such that
\begin{equation}
W_{1}\left(x\right)\leq V\left(x\left(t\right),t\right)\leq W_{2}\left(x\right)\:\:\:\:\forall t\geq0,\:\:\forall x\in\mathcal{D},\label{eq: v_bound}
\end{equation}
\begin{equation}
\dot{V}\left(x\left(t\right),t\right)\overset{a.e.}{\leq}-W\left(x\right)\label{eq: v_tilde_dot}
\end{equation}
where $W_{1}\left(x\right)$ and $W_{2}\left(x\right)$ are continuous
positive definite functions, $W\left(x\right)$ is a continuous positive
semi-definite function on $\mathcal{D}$, choose $r>0$ and $c>0$
such that $B_{r}\subset\mathcal{D}$ and $c<\min\limits _{\left\Vert x\right\Vert =r}W_{1}\left(x\right)$
and $x\left(t\right)$ is a Filippov solution to (\ref{eq: dynamics})
where $x\left(t_{0}\right)\in\left\{ x\in B_{r}\:|\: W_{2}\left(x\right)\leq c\right\} $.
Then $x\left(t\right)$ is bounded and satisfies
\begin{equation}
W\left(x\left(t\right)\right)\rightarrow0\:\:\:\: as\:\:\:\: t\rightarrow\infty.\label{eq: w zero 1}
\end{equation}
\end{cor}
\begin{IEEEproof}
Since $x\left(t\right)$ is a Filippov solution to (\ref{eq: dynamics}),
$\left\{ x\in B_{r}\:|\: W_{1}\left(x\right)\leq c\right\} $ is in
the interior of $B_{r}$. Define a time-dependent set $\Omega_{t,c}$
by
\[
\Omega_{t,c}=\left\{ x\in B_{r}\:|\: V\left(x,t\right)\leq c\right\} .
\]
From (\ref{eq: v_bound}), the set $\Omega_{t,c}$ contains $\left\{ x\in B_{r}\:|\: W_{2}\left(x\right)\leq c\right\} $
since
\[
W_{2}\left(x\right)\leq c\Rightarrow V\left(x,t\right)\leq c.
\]
On the other hand, $\Omega_{t,c}$ is a subset of $\left\{ x\in B_{r}\:|\: W_{1}\left(x\right)\leq c\right\} $
since
\[
V\left(x,t\right)\leq c\Rightarrow W_{1}\left(x\right)\leq c.
\]
Thus, 
\[
\left\{ x\in B_{r}\:|\: W_{2}\left(x\right)\leq c\right\} \subset\Omega_{t,c},
\]
\[
\Omega_{t,c}\subset\left\{ x\in B_{r}\:|\: W_{1}\left(x\right)\leq c\right\} \subset B_{r}\subset\mathcal{D}.
\]

Based on (\ref{eq: v_tilde_dot}),\textcolor{red}{{} }$\dot{V}\left(x(t),t\right)\overset{a.e.}{\leq}0$,
hence, $V\left(x\left(t\right),t\right)$ is non-increasing from Lemma
\ref{lem: v stable}. For any $t_{0}\geq0$ and any $x\left(t_{0}\right)\in\Omega_{t_{0},c}$,
the solution starting at $\left(x\left(t_{0}\right),t_{0}\right)$
stays in $\Omega_{t,c}$ for every $t\geq t_{0}$. Therefore, any
solution starting in $\left\{ x\in B_{r}\:|\: W_{2}\left(x\right)\leq c\right\} $
stays in $\Omega_{t,c}$, and consequently in $\left\{ x\in B_{r}\:|\: W_{1}\left(x\right)\leq c\right\} $,
for all future time. Hence, the Filippov solution $x\left(t\right)$
is bounded such that $\left\Vert x\left(t\right)\right\Vert <r$,
$\forall t\geq t_{0}$. 

From Lemma \ref{lem: v stable}, $V\left(x\left(t\right),t\right)$
is also bounded such that $V\left(x\left(t\right),t\right)\leq V\left(x\left(t_{0}\right),t_{0}\right)$.
Since $\dot{V}\left(x(t),t\right)$ is Lebesgue measurable from (\ref{eq: v_tilde_dot}),
\begin{equation}
\int_{t_{0}}^{t}W\left(x\left(\tau\right)\right)d\tau\leq-\int_{t_{0}}^{t}\dot{V}\left(x\left(\tau\right),\tau\right)d\tau,\label{eq: int_W}
\end{equation}
\[
-\!\!\int_{t_{0}}^{t}\!\!\dot{V}\!\!\left(x\!\left(\tau\right)\!,\tau\right)\! d\tau\!=\! V\!\!\left(x\!\left(t_{0}\right)\!,t_{0}\right)\!-\! V\!\!\left(x\!\left(t\right)\!,t\right)\!\leq\! V\!\!\left(x\!\left(t_{0}\right)\!,t_{0}\right).
\]
Therefore, $\int_{t_{0}}^{t}W\left(x\left(\tau\right)\right)d\tau$
is bounded $\forall t>t_{0}$. Existence of $\lim\limits _{t\rightarrow\infty}\int_{t_{0}}^{t}W\left(x\left(\tau\right)\right)d\tau$
is guaranteed since the left-hand side of (\ref{eq: int_W}) is monotonically
nondecreasing (based on the definition of $W\left(x\right)$) and
bounded above. Since every absolutely continuous function is uniformly
continuous, $x\left(t\right)$ is uniformly continuous. Because $W\left(x\right)$
is continuous in $x$, and $x$ is on the compact set $B_{r}$, $W\left(x\left(t\right)\right)$
is uniformly continuous in $t$ on $\left(t_{0},\infty\right]$. Therefore,
by Lemma \ref{lem: barbalats}, 
\begin{equation}
W\left(x\left(t\right)\right)\rightarrow0\: as\: t\rightarrow\infty.\label{eq: w zero}
\end{equation}
\end{IEEEproof}
\begin{rem}
From Def. \ref{def: filippov}, $K\left[f\right]\left(x,t\right)$
is an upper semi-continuous, nonempty, compact and convex valued map.
While existence of a Filippov solution for any arbitrary initial condition
$x\left(t_{0}\right)\in\mathcal{D}$ is provided by the definition,
generally speaking, the solution is non-unique\cite{Filippov1988,Aubin1984}. 

Note that Corollary \ref{thm: main_result} establishes (\ref{eq: w zero})
for a specific $x\left(t\right)$. Under the stronger condition that%
\footnote{The inequality $\dot{\tilde{V}}\left(x,t\right)\leq W\left(x\right)$
is used to indicate that every element of the set $\dot{\tilde{V}}\left(x,t\right)\mbox{ is less than or equal to the scalar \ensuremath{W\left(x\right)}. }$%
} $\dot{\tilde{V}}\left(x,t\right)\leq W\left(x\right)\:\forall x\in\mathcal{D}$,
it is possible to show that (\ref{eq: w zero}) holds for all Filippov
solutions of (\ref{eq: dynamics}). The next corollary is presented
to illustrate this point.\end{rem}
\begin{cor}
\label{thm: main_result_2} For the system given in (\ref{eq: dynamics}),
let $\mathcal{D}\subset\mathbb{R}^{n}$ be a domain containing $x=0$
and suppose $f$ is Lebesgue measurable and essentially locally bounded
on $\mathcal{D}\times\left[0,\infty\right)$. Furthermore, suppose
$f\left(0,t\right)$ is uniformly bounded for all $t\geq0$. Let $V:\mathcal{D}\times\left[0,\infty\right)\rightarrow\mathbb{R}$
be continuously differentiable in x, locally Lipschitz in t, and regular
such that
\begin{equation}
W_{1}\left(x\right)\leq V\left(x,t\right)\leq W_{2}\left(x\right)\label{eq: v_bound_2}
\end{equation}
\begin{equation}
\dot{\tilde{V}}\left(x,t\right)\leq-W\left(x\right)\label{eq: v_tilde_dot_2}
\end{equation}
$\forall t\geq0,\:\:\forall x\in\mathcal{D}$ where $W_{1}\left(x\right)$
and $W_{2}\left(x\right)$ are continuous positive definite functions,
and $W\left(x\right)$ is a continuous positive semi-definite function
on $\mathcal{D}$. Choose $r>0$ and $c>0$ such that $B_{r}\subset\mathcal{D}$
and $c<\min\limits _{\left\Vert x\right\Vert =r}W_{1}\left(x\right)$.
Then, all Filippov solutions of (\ref{eq: dynamics}) such that $x\left(t_{0}\right)\in\left\{ x\in B_{r}\:|\: W_{2}\left(x\right)\leq c\right\} $
are bounded and satisfy
\begin{equation}
W\left(x\left(t\right)\right)\rightarrow0\:\: as\:\: t\rightarrow\infty.\label{eq: w zero 2}
\end{equation}
\end{cor}
\begin{IEEEproof}
Let $x\left(t\right)$ be any arbitrary Filippov solution of (\ref{eq: dynamics}).
Then, from Lemma \ref{lem: chain rule}, and (\ref{eq: v_tilde_dot_2}),
$\dot{V}\left(x\left(t\right),t\right)\overset{a.e.}{\leq}-W\left(x\left(t\right)\right)$,
which is the condition (\ref{eq: v_tilde_dot}). Since the selection
of $x\left(t\right)$ is arbitrary, Corollary \ref{thm: main_result}
can be used to imply that the result in (\ref{eq: w zero}) holds
for each $x\left(t\right)$. Hence, Corollary \ref{thm: main_result_2}
holds.\end{IEEEproof}
\begin{rem}
In the case of some systems (e.g., closed loop error systems with
sliding mode control laws), it may be possible to show that Corollary
\ref{thm: main_result_2} is more easily applied. However, in other
cases, it may be difficult to satisfy the inequality in (\ref{eq: v_tilde_dot_2}).
The usefulness of Corollary \ref{thm: main_result} is demonstrated
in those cases where it is difficult or impossible to show that the
inequality in (\ref{eq: v_tilde_dot_2}) can be satisfied, but it
is possible to show that (\ref{eq: v_tilde_dot}) can be satisfied
for almost all time.
\end{rem}

\section{Conclusion}

In this note, the Lasalle-Yoshizawa Theorem is extended to differential
systems whose right-hand sides are discontinuous in the state and
piecewise continuous in time. The result presents two theoretical
tools applicable to nonautonomous systems with discontinuities in
the closed-loop error system. Generalized Lyapunov-based analysis
methods are developed utilizing differential inclusions in the sense
of Filippov to achieve asymptotic convergence when the candidate Lyapunov
derivative is upper bounded by a negative semi-definite function.

\bibliographystyle{IEEEtran}
\bibliography{ncr,master}

\end{document}